\documentclass[10pt,twoside]{amsart}
\usepackage{amsmath}
\usepackage{amssymb}
\def\firstpage{1}\def\lastpage{10}

\makeatletter
\def\magnification{\afterassignment\m@g\count@}
\def\m@g{\mag=\count@\hsize6.5truein\vsize8.9truein\dimen\footins8truein}
\makeatother

\font\eightrm=cmr8
\font\caps=cmcsc10                    % Theorem, Lemma etc
\font\scaps=cmcsc8

\font\Caps=cmcsc10 scaled \magstep1   % Title

\makeatletter
\@specialpagefalse\headheight=8.5pt
\def\DocMath{{\def\th{\thinspace}\scaps Doc.\th Math.\th J.\th DMV}}
\renewcommand{\@oddfoot}{\hfill\scaps Documenta Mathematica $\cdot$ Extra
        Volume ICM  1998 $\cdot$ \number\firstpage--\lastpage\hfill}
\renewcommand{\@evenfoot}{\ifnum\thepage>\lastpage\hfill\scaps
    Documenta Mathematica $\cdot$ Extra Volume ICM  1998 $\cdot$
\hfill\else\@oddfoot\fi}%
\renewcommand{\@evenhead}{%
    \ifnum\thepage>\lastpage\rlap{\thepage}\hfill%
    \else\rlap{\thepage}\slshape\leftmark\hfill\caps\SAuthor\hfill\fi}%
\renewcommand{\@oddhead}{%
    \ifnum\thepage=\firstpage{\DocMath\hfill\llap{\thepage}}%
    \else{\slshape\rightmark}\hfill\caps\STitle\hfill\llap{\thepage}\fi}%
\makeatother

\def\TSkip{\bigskip}
\newbox\TheTitle{\obeylines\gdef\GetTitle #1
\ShortTitle  #2
\SubTitle    #3
\Author      #4
\ShortAuthor #5
\EndTitle
{\setbox\TheTitle=\vbox{\baselineskip=20pt\let\par=\cr\obeylines%
\halign{\centerline{\Caps##}\cr\noalign{\medskip}\cr#1\cr}}%
        \copy\TheTitle\TSkip\TSkip%
\def\next{#2}\ifx\next\empty\gdef\STitle{#1}\else\gdef\STitle{#2}\fi%
\def\next{#3}\ifx\next\empty%
    \else\setbox\TheTitle=\vbox{\baselineskip=20pt\let\par=\cr\obeylines%
    \halign{\centerline{\caps##} #3\cr}}\copy\TheTitle\TSkip\TSkip\fi%
\centerline{\caps #4}\TSkip\TSkip%
\def\next{#5}\ifx\next\empty\gdef\SAuthor{#4}\else\gdef\SAuthor{#5}\fi%
\catcode'015=5}}\def\Title{\obeylines\GetTitle}
\def\Abstract{\begingroup\narrower
    \parskip=\medskipamount\parindent=0pt{\caps Abstract. }}
\def\EndAbstract{\par\endgroup\TSkip}

\long\def\MSC#1\EndMSC{\def\arg{#1}\ifx\arg\empty\relax\else
     {\par\narrower\noindent%
     1991 Mathematics Subject Classification: #1\par}\fi}

\long\def\KEY#1\EndKEY{\def\arg{#1}\ifx\arg\empty\relax\else
        {\par\narrower\noindent Keywords and Phrases: #1\par}\fi\TSkip}

\newbox\TheAdd\def\Addresses{\vfill\copy\TheAdd\vfill
    \ifodd\number\lastpage\vfill\eject\phantom{.}\vfill\eject\fi}
{\obeylines\gdef\GetAddress #1
\Address #2
\Address #3
\Address #4
\EndAddress
{\def\xs{6truecm}
\setbox0=\vtop{{\obeylines\hsize=\xs#1\par}}\def\next{#2}
\ifx\next\empty % 1 address
     \setbox\TheAdd=\hbox to\hsize{\hfill\copy0\hfill}
\else\setbox1=\vtop{{\obeylines\hsize=\xs#2\par}}\def\next{#3}
\ifx\next\empty % 2 addresses
     \setbox\TheAdd=\hbox to\hsize{\hfill\copy0\hfill\copy1\hfill}
\else\setbox2=\vtop{{\obeylines\hsize=\xs#3\par}}\def\next{#4}
\ifx\next\empty\ % 3 addresses
     \setbox\TheAdd=\vtop{\hbox to\hsize{\hfill\copy0\hfill\copy1\hfill}
                \vskip20pt\hbox to\hsize{\hfill\copy2\hfill}}
\else\setbox3=\vtop{{\obeylines\hsize=\xs#4\par}}
     \setbox\TheAdd=\vtop{\hbox to\hsize{\hfill\copy0\hfill\copy1\hfill}
                \vskip20pt\hbox to\hsize{\hfill\copy2\hfill\copy3\hfill}}
\fi\fi\fi\catcode'015=5}}\gdef\Address{\obeylines\GetAddress}

\begin{document}
\theoremstyle{plain}
\newtheorem{thm}{Theorem}[section]
\newtheorem{cor}[thm]{Corollary}
\newtheorem{conj}{Conjecture}
\renewcommand{\theconj}{}

\def \CPb {\overline{\mathbf{CP}}{}^{\hspace{.005in} 2}}
\def \CP {{\mathbf{CP}}^{\hspace{.005in} 2}}
\def \la {\langle}
\def \ra {\rangle}
\def \Z {\mathbf{Z}}
\def \R {\mathbf{R}}
\def \Sig{\Sigma}
\def \vt {\vartheta}
\def \a {\alpha}
\def \b {\beta}
\def \o {\omega}
\def \s {\sigma}
\def \t {\tau}
\def \bd {\partial}
\def \x {\times}
\def \ve {\varepsilon}
\def \e {\epsilon}
\def \ssw {\text{SW}}
\def \sw {\mathcal{SW}}
\def \al {\text{A}}
\def \DD {\Delta}
\def \DN {\nabla}
\def \ggr {\text{Gr}}
\def \gr {\mathcal{G}r}
%%%%% ------------- fill in your data below this line  -------------------
%%%%%    The following lines \Title ... \EndAddress must ALL be present
%%%%%    and in the given order.
\Title
     Constructions of smooth 4-manifolds
%%%%%    Put here the title. Line breaks will be recognized.
\ShortTitle
%%%%%    Running title for odd numbered pages, ONE line, please.
%%%%%    If none is given, \Title will be used instead.
\SubTitle
%%%%%    A possible subtitle goes here.
\Author
         Ronald Fintushel%
      \footnote{\eightrm Partially supported by NSF Grant DMS9704927}\
        and Ronald J. Stern%
      \footnote{\eightrm Partially supported by NSF Grant DMS9626330}
%%%%%    Put here name(s) of authors. Line breaks will be recognized.
\ShortAuthor R. Fintushel and R. J. Stern
%%%%%%   Running title for even numbered pages, ONE line, please.
%%%%%%   If none is given, \Author will be used instead.
\EndTitle
\Abstract
%%%%%    Put here the abstract of your manuscript.
We describe a collection of
constructions which illustrate a panoply of ``exotic'' smooth 4-manifolds.
\EndAbstract
\MSC
%%%%%    1991 Mathematics Subject Classification:
     57R55
\EndMSC
\KEY
%%%%%    Keywords and Phrases:
\EndKEY
%%%%%    All 4 \Address lines below must be present. To center the last
%%%%%    entry, no empty lines must be between the following \Address
%%%%%    and \EndAddress lines.
\Address
               Ronald Fintushel
               Department of Mathematics
               Michigan State University
               East Lansing, Michigan 48824
               U.S.A.
               ronfint@math.msu.edu
%%%%%    Address of first Author here
\Address       Ronald J. Stern
               Department of Mathematics
               University of California
               Irvine, California 92697
               U.S.A.
               rstern@uci.edu
%%%%%    Address of second Author here etc.
\Address
\Address
\EndAddress
%%
%%       Make sure the last tex command in your manuscript
%%       before the first \end{document} is the command  \Addresses
%%
%%---------------------Here the prologue ends---------------------------------
%%--------------------Here the manuscript starts------------------------------

\section{Introduction\label{Intro}}

At the time of the previous (1994) International Congress of
Mathematicians, steady, but slow, progress was being made on the
classification of simply connected closed smooth 4-manifolds. In
particular, the Donaldson invariants had begun to take a particularly nice
form \cite{KM} (also \cite{FSstructure}), their computations were
becoming more routine \cite{rat}, and their behavior under blowing up
(i.e. taking connected sum with $\CPb$)
 was well understood \cite{blowup}. Due to the complexity of the Donaldson
invariants, great hope was held out that an even better understanding of
these invariants would close the books on the classification of
simply connected 4-manifolds.

A few short months after the 1994 ICM, the 4-manifold community was
blind-sided by the introduction of the now famous Seiberg-Witten equations
\cite{W}. Most of the results obtained by using  Donaldson theory were
found to
have quicker, and sometimes more general, counterparts using the Seiberg-Witten
technology. The potential applications of the difficult Donaldson
technology became much more transparent using these new equations. As of
July 1998, there is good news as well as bad news. The good news is that
many of the earlier focus problems have been solved. In particular, the
Thom conjecture \cite{KM2} and  its natural generalizations have been
verified \cite{MST,OS}; also the study of symplectic 4-manifolds
has taken a more central role \cite{T1,T2,TSWG,TCPH}.
The bad news is that recent constructions and computations indicate that
the Seiberg-Witten and Donaldson theories are too weak to distinguish
simply connected smooth 4-manifolds \cite{KL4M}.  It is these
latter constructions and computations that we will discuss at this 1998
International Congress of Mathematicians. It is becoming more
apparent that we are seeing only a small constellation of 4-dimensional
manifolds. More seriously, we are lacking a reasonable conjectural
classification of simply connected closed smooth 4-manifolds.

Current technology has given us many more
4-manifolds than had been expected in 1994. The authors hope that during
the 2002 ICM the
construction of large classes of new 4-manifolds will be discussed; in
particular, they hope that a sufficiently large collection of 4-manifolds
will have been
discovered so as to allow for
some general patterns to emerge and, at least, a conjectural classification
to again be on the books.

\section{The knot surgery construction \label{knots}}

Let $X$ be a simply connected oriented smooth closed 4-manifold. Its most basic
invariant
is its intersection form
\[ Q_X:H_2(X;\Z)\otimes H_2(X;\Z)\to\Z\]
defined by counting signed transverse intersections of embedded oriented
surfaces
representing given homology classes. It is a famous theorem of M. Freedman
\cite{Fr} that $Q_X$ determines the homeomorphism type of $X$, and an equally
renowned theorem of S.K. Donaldson \cite{Donpoly} that $Q_X$ is not
sufficient to
determine the diffeomorphism type of $X$. In this section we shall discuss
geometric operations on a given smooth 4-manifold which preserve the underlying
topological structure and alter its smooth structure. In particular, we shall
consider the following construction: Let $X$ be a simply connected smooth
4-manifold which contains a smoothly embedded torus
$T$ of self-intersection 0. Given a knot $K$ in $S^3$, we replace a
tubular neighborhood of $T$ with $S^1\times (S^3\setminus K)$ to
obtain the {\it knot surgery manifold} $X_K$.

More formally, this procedure is accomplished by performing 0-framed
surgery on $K$ to obtain the 3-manifold $M_K$. The meridian $m$ of $K$ can be
viewed as a circle in $M_K$; so in $S^1\times M_K $ we have the smooth torus
$T_m=S^1\times m$ of self-intersection~0. Since a neighborhood of
$m$ has a canonical framing in $M_K$, a neighborhood of the torus
$T_m$ in $S^1\times M_K$ has a canonical identification with
$T_m\times D^2$. The knot surgery manifold $X_K$ is given by the
fiber sum \[ X_K=X\#_{T=T_m}S^1\times M_K=(X\setminus T\times D^2)
\cup(S^1\times M_K\setminus T_m\times D^2)\] where the two pieces are
glued together so as to preserve the homology class
$[{\rm{pt}}\times \partial D^2]$. This latter condition does not, in
general, completely determine the isotopy type of the gluing, and
$X_K$ is taken to be any manifold constructed in this fashion.

Because $S^1\times (S^3\setminus K)$ has the same homology as a tubular
neighborhood of $T$ in $X$ (and because the gluing preserves
$[{\rm{pt}}\times \partial D^2]$)
the homology and intersection form of
$X_K$ will agree
with that of $X$. If it is also assumed that $X\setminus T$ is simply
connected,
then $\pi_1(X_K)=1$; so $X_K$ will be homeomorphic to $X$.

In order to distinguish the diffeomorphism types of the $X_K$, we rely on
Seiberg-Witten invariants. We view the Seiberg-Witten invariant of a smooth
4-manifold as a multivariable (Laurent) polynomial. To do this, recall
that the Seiberg-Witten
invariant  of a smooth closed oriented $4$-manifold
$X$ with $b_2 ^+(X)>1$ is an integer-valued function which is defined on
the set of $spin ^{\, c}$ structures over $X$ (cf. \cite{W}). In case
$H_1(X,\Z)$
has no 2-torsion (for example, as here where $X$ is simply connected) there
is a
natural identification of the $spin ^{\, c}$ structures of
$X$ with the characteristic elements of $H_2(X,\Z)$ (i.e. those elements
$k$ whose
Poincar\'e duals $\hat{k}$ reduce mod~2 to $w_2(X)$).
In this case we view the
Seiberg-Witten invariant as
\[ \ssw_X: \lbrace k\in H_2(X,\Z)|\hat{k}\equiv w_2(TX)\pmod2)\rbrace
\rightarrow \Z. \]
The sign of $\sw_X$
depends on an orientation of
$H^0(X,\R)\otimes\det H_+^2(X,\R)\otimes \det H^1(X,\R)$. If $\ssw_X(\b)\neq
0$, then $\b$ is called a {\it{basic class}} of $X$. It is a fundamental
fact that the set of
basic classes is
finite. Furthermore, if $\b$ is a basic class, then so is $-\b$ with
$\ssw_X(-\b)=(-1)^{(\text{e}+\text{sign})(X)/4}\,\ssw_X(\b)$ where
$\text{e}(X)$ is
the Euler number and $\text{sign}(X)$ is the signature of $X$.

Now let
$\{\pm \b_1,\dots,\pm \b_n\}$ be the set of nonzero basic classes for $X$.
Consider variables $t_{\b}=\exp(\b)$ for each $\b\in H^2(X;\Z)$ which
satisfy the
relations $t_{\a+\b}=t_{\a}t_{\b}$. We may then view the Seiberg-Witten
invariant
of $X$ as the Laurent polynomial
\[\sw_X = \ssw_X(0)+\sum_{j=1}^n \ssw_X(\b_j)\cdot
(t_{\b_j}+(-1)^{(\text{e}+\text{sign})(X)/4}\, t_{\b_j}^{-1}).\]
As an example of this notational device, consider the simply connected
minimally elliptic
surface
$E(n)$ with holomorphic Euler characteristic $n$ and no multiple fibers. Its
Seiberg-Witten invariant is
$\sw_{E(n)}=(t-t^{-1})^{n-2}$ where $t=t_F$ for $F$ the fiber class. Thus,
$\ssw_{E(n)}((n-2m)F)=(-1)^{m-1}\binom{n-2}{m-1}$ for
$m=1,\dots,n-1$ and $\ssw_{E(n)}(\a)=0$ for any other $\a$.
When $b^+(X)>1$, the Laurent polynomial $\sw_X$ is a diffeomorphism
invariant of $X$.

For our theorem, we need to place a mild hypothesis on the embedded
torus $T$. We say that a smoothly embedded torus representing a
nontrivial homology class $[T]$ is {\it{c-embedded}} if there is a
neighborhood $N$ of $T$ in $X$ and a diffeomorphism $\varphi:N\to U$ where
$U$ is a neighborhood of a cusp fiber in an elliptic surface and
$\varphi(T)$ is a
smooth elliptic fiber in $U$. Equivalently, $T$ is c-embedded if it
contains two
simple closed curves which generate $\pi_1(T)$ and which bound  vanishing
cycles
in $X$. Note that a c-embedded torus has self-intersection~0.

\begin{thm}[\cite{KL4M}]\label{knotthm} Let $X$ be a simply connected
oriented smooth
4-manifold with $b^+>1$. Suppose that $X$ contains a c-embedded
torus $T$ with $\pi_1(X\setminus T)=1$, and let $K$ be any knot in $S^3$.
Then the knot
surgery manifold $X_K$ is homeomorphic to $X$ and has Seiberg-Witten invariant
\[ \sw_{X_K}=\sw_{X}\cdot\DD_K(t) \]
where $\DD_K(t)$ is the symmetrized Alexander polynomial of $K$ and
$t=\exp(2[T])$.
\end{thm}

For example, the theorem applies to the K3-surface $E(2)$ where $T$ is a
smooth elliptic
fiber, and since $\sw_{E(2)}=1$, we have $\sw_{E(2)_K}=\DD_K(t)$. It is a
theorem of
Seifert that any Laurent polynomial of the
form $P(t)=a_0+\sum\limits_{j=1}^na_j(t^j+t^{-j})$ with coefficient sum
$P(1)=\pm1$ is the Alexander polynomial of some
knot in $S^3$.
Call such a Laurent polynomial an {\it {$A$-polynomial}}. It follows that
if $(X,T)$
satisfies the hypothesis of Theorem~\ref{knotthm}, then for any
$A$-polynomial $P(t)$, there is a smooth simply connected 4-manifold
$X_{P}$ which is
homeomorphic to $X$
and has Seiberg-Witten invariant $\sw_{X_{P}}=\sw_X\cdot P(t)$
where $t=\exp(2[T])$. In particular, for each $A$-polynomial $P(t)$, there
is a manifold
homeomorphic to the K3-surface with $\sw=P(t)$.

The relationship between Seiberg-Witten type invariants and the Alexander
polynomial was
first discovered by Meng and Taubes. In \cite{MT} they showed that the
3-manifold
Seiberg-Witten invariant is related to Milnor torsion.

If one starts with a fibered knot $K$, then $S^1\x M_K$ is a surface bundle
over a torus and thus carries a symplectic structure \cite{Th} for which
$T_m$ is a
symplectic submanifold. Thus if
$X$ is a symplectic 4-manifold
containing a c-embedded symplectic torus $T$, then
$X_K=X\#_{T=T_m}S^1\times M_K$ is also
symplectic \cite{Gompf, MW}. In a fashion similar
to the treatment of the Seiberg-Witten invariant as a Laurent polynomial,
one can
view the Gromov invariant of a symplectic 4-manifold $X$ as a polynomial
$\gr_X = \sum\ggr_X(\b)\,t_{\b}$
where $\ggr_X(\b)$ is the usual Gromov invariant of $\b$. Let
$\al_K(t)=t^{d}\DD_K(t)$
denote the normalized Alexander polynomial, where $d$ is the degree of
$\DD_K(t)$.
As a corollary to Theorem~\ref{knotthm} and the theorems of Taubes relating the
Seiberg-Witten and Gromov invariants of a symplectic
$4-$manifold \cite{TSWG,TCPH} we have:

\begin{cor}[\cite{KL4M}] \label{symcor} Let $X$ be a symplectic
$4$-manifold with
$b^+>1$ containing a symplectic c-embedded torus $T$. If
$K$ is a fibered knot, then $X_K$ is a symplectic $4$-manifold whose Gromov
invariant is \
$\gr_{X_K}=\gr_X\cdot\al_K(\t)$ \
where $\t=\exp([T])$.
\end{cor}

This last calculation can also be made purely within the realm of
symplectic topology
\cite{IP,L}. Our interest is directed more to the opposite situation. The
Alexander
polynomial of a fibered knot is monic; i.e. its top coefficient is $\pm1$.
On the other
hand:

\begin{cor}[\cite{KL4M}] \label{nonsymp} If $\DD_K(t)$ is not monic, then
$X_K$ does
not admit a symplectic structure. Furthermore, if
$X$  contains a homologically nontrivial surface $\Sigma_g$ of genus $g$
disjoint from
$T$  with $[\Sigma_g]^2 < 2-2g$ if $g>0$ or $[\Sigma_g]^2 < 0$ if
$g=0$, then $X_K$ with the opposite orientation does not admit a symplectic
structure.
\end{cor}

Until the summer of 1996, it was still a plausible conjecture (sometimes
called the
`minimal conjecture') that each irreducible simply connected 4-manifold
should admit a
symplectic structure with one of its orientations. The first
counterexamples to this
conjecture were constructed by Z. Szabo \cite{S}.
The knot surgery construction gives a multitude of examples of simply
connected irreducible
`nonsymplectic' 4-manifolds. In fact, if $X$ is simply connected with
$\sw_X\ne0$ and if
$X$ contains a c-embedded torus $T$ with $\pi_1(X\setminus T)=1$, then
Theorem~\ref{knotthm} and  Corollary~\ref{nonsymp} imply that there are
infinitely many
distinct nonsymplectic smooth 4-manifolds $X_K$ homeomorphic to $X$.

If $K_1$ and $K_2$ have the same Alexander polynomial, Seiberg-Witten
invariants are
not able to distinguish $X_{K_1}$ from $X_{K_2}$. For example, take
$X=E(2)$. Then $X_K$ has
a self-intersection 0 homology class $\s$ satisfying $\s \cdot [T]=1$ which is
represented by an embedded surface of genus $g(K)+1$ where $g(K)$ is the
genus of $K$. One
might hope that these classes could be used to distinguish $X_{K_1}$ from
$X_{K_2}$ when
$g(K_1)\ne g(K_2)$.

\begin{conj} For $X=E(2)$, the manifolds  $X_{K_1}$ and $X_{K_2}$ are
diffeomorphic if
and only if $K_1$ and $K_2$ are equivalent knots.
\end{conj}

The proof of Theorem~\ref{knotthm} proceeds by successively simplifying the
manifold $X_K$
in a fashion which mimics the calculation of the Alexander polynomial of $K$
via skein
relations. Recall that $\DD_K(t)$ can be calculated via the relation
\begin{equation} \DD_{K_+}(t)=\DD_{K_-}(t)+(t^{1/2}-t^{-1/2})\cdot\DD_{K_0}(t)
\label{boogie}\end{equation}
where $K_+$ is an oriented knot or link, $K_-$ is the result of
changing a single oriented positive (right-handed) crossing in $K_+$ to a
negative (left-handed) crossing, and $K_0$ is the result of resolving the
crossing
as shown in Figure~1.

\centerline{\unitlength .65cm
\begin{picture}(9,4)
\put (1,1){\vector(1,2){1.12}}
\put (2,1){\line(-1,2){0.45}}
\put (1.45,2.1){\vector(-1,2){.575}}
\put (1.25,0.5){$K_+$}
\put (4,1){\line(1,2){0.45}}
\put (4.55,2.1){\vector(1,2){.575}}
\put (5,1){\vector(-1,2){1.12}}
\put (4.25,0.5){$K_-$}
\put (7,1){\line(1,4){.255}}
\put (8,1){\line(-1,4){.255}}
\put (7.255,2){\vector(-1,4){.3}}
\put (7.745,2){\vector(1,4){.3}}
\put (7.25,0.5){$K_0$}
\put (3.3,-.15){Figure 1}
\end{picture}}
\medskip

The point of using \eqref{boogie} to calculate $\DD_K$ is that $K$ can be
simplified to an unknot via a sequence of crossing changes. One builds a
`resolution tree' starting from $K$ and at each stage adding the
bifurcation of Figure~2,
where each $K_+$, $K_-$, $K_0$ is a knot or 2-component link, and
so that at the bottom of the tree, there are only unknots, and split links.
Then,
because the Alexander polynomial of an unknot is $1$, and is $0$ for a
split link
(of more than one component) one can can work backwards using
\eqref{boogie} to calculate
$\DD_K(t)$.

\centerline{\unitlength .5cm
\begin{picture}(5.5,5.5)
\put (2.25,4){\line(-3,-4){1.5}}
\put (2.75,4){\line(3,-4){1.5}}
\put (2.25,4.3){$K_+$}
\put (0.25,1.2){$K_-$}
\put (4,1.2){$K_0$}
\put (1.4,.3){Figure 2}
\end{picture}}

The manifold $X_{K_+}$ can be obtained from $X_{K_-}$ by means of a
$(+1)$-log transform on
a nullhomologous torus in $X_{K_-}$, and then the gluing theorems of
\cite{MMS} show that
$\sw_{X_{K_+}}$ can be computed in terms of the Seiberg-Witten invariants
of $X_{K_-}$  and
a manifold $X_{K_-,0}$ obtained by a  0-log transform on $X_{K_-}$. With
some work, this
leads to a related resolution diagram of 4-manifolds where each knot $K'$
corresponds
to $X_{K'}$, and this diagram can be used to prove Theorem~\ref{knotthm}.

We conclude this section by pointing out that the knot surgery construction
can be
generalized to manifolds with $b^+=1$ and to links in $S^3$ of more than one
component in a
more-or-less obvious way. One glues the complements of c-embedded tori in
4-manifolds to
the product of $S^1$ with the link complement. See \cite{KL4M} for details.
For example, if
to each boundary component of $S^1\x (S^3\setminus N(L))$ we glue $E(1)$
minus the
neighborhood of a smooth elliptic fiber, we obtain a manifold with
$\sw=\DD_L(t_1,\dots,t_n)$, the multivariable Alexander polynomial of the
link. Szabo's
examples in \cite{S} can be obtained from this construction.

\section{Embeddings of surfaces in 4-manifolds \label{surf}}

Knot surgery can also be used to change the embedding of a surface in a
fixed 4-manifold.
To motivate the construction, note that one can tie a knot in the core $\{
0\}\x I$ of a
cylinder $D^2\x I$ by removing a tubular neighborhood of a meridian circle
and replacing it
with a knot complement $S^3\setminus K$. We shall perform a parametrized
version of this
construction in the 4-manifold setting.  Consider an oriented surface $\Sig$ of
genus $g>0$ which is smoothly embedded in a simply connected 4-manifold
$X$. Let $\a$ be a
simple closed curve on $\Sig$ which is part of a symplectic basis, and let
$\a\x I$ be an
annular neighborhood of $\a$ in $\Sig$. In $X$ we see the neighborhood
$D^2\x\a\x I$. For a
fixed knot $K$ in $S^3$, we parametrize the above construction so as to
perform it
on each of the cylinders
$D^2\x\{ y\}\x I$, $y\in\a$, to obtain an embedded surface $\Sig_K$. This
is equivalent to
performing knot surgery on the (nullhomologous) rim torus $R=\bd D^2\x\a$.
We call this
operation {\it{rim surgery}}.

\begin{thm}[\cite{surfaces}]\label{homeo} Let $X$ be a simply connected
smooth 4-manifold
with an embedded surface $\Sig$ of positive genus. Suppose that
$\pi_1(X\setminus\Sig)=1$. Then for each knot
$K$ in
$S^3$, rim surgery produces a surface $\Sig_K$, and there is a
homeomorphism $(X, \Sig)\cong (X,\Sig_K)$.
\end{thm}

The Seiberg-Witten invariant can be used to study these embeddings, but
first, an
auxilliary construction is needed. For each positive integer $g$, let $Y_g$ be
the union of the Milnor fiber of the $(2,2g+1,4g+1)$ Brieskorn singularity
and a
generalized nucleus consisting of the 4-manifold obtained as the trace of
the 0-framed
surgery on $(2,2g+1)$ torus knot in $\bd B^4$ and a $-1$ surgery on a
meridian. Then $Y_g$
is a Kahler surface and admits a holomorphic fibration over
${\mathbf{CP}}^{1}$ with
generic fiber a surface $S_g$ of genus $g$.

Let $(X,\Sig)$ be as in Theorem~\ref{homeo}, and suppose that the
self-intersection $\Sig^2=0$.  We call $(X,\Sig)$ an {\it SW-pair}\/ if
satisfies the property that $\sw_{X\#_{\Sig=S_g} Y_g}\ne 0$. (In general, if
$\Sig^2=n>0$, one makes this definition by first blowing up $n$ times.)
For example, if $X$ is symplectic and
$\Sig$ is a
symplectic submanifold (of square 0), then $X\#_{\Sig=S_g} Y_g$ is
symplectic, and it
follows that
$(X,\Sig)$ is an SW-pair. In $X\#_{\Sig=S_g} Y_g$, the rim torus $R$
becomes homologically
essential and is c-embedded. We can use Theorem~\ref{knotthm} to calculate
Seiberg-Witten
invariants:
\[\sw_{X\#_{\Sig_K=S_g}Y_g}=\sw_{(X\#_{\Sig=S_g}Y_g)_K}=
\sw_{X\#_{\Sig=S_g}Y_g}\cdot\DD_K(r)\]
where $r=\exp(2[R])$, viewing $[R]$ as a class in the fiber sum. We have:
\begin{thm}[\cite{surfaces}]\label{ss} Consider any SW-pair $(X,\Sig)$ with
$\Sig^2\ge0$. If
$K_1$ and $K_2$ are two knots in $S^3$ and if there is a diffeomorphism of
pairs
$(X,\Sig_{K_1})\cong (X,\Sig_{K_2})$, then  $\DD_{K_1}(t)=\DD_{K_2}(t)$.
\end{thm}

As a special case:
\begin{thm}[\cite{surfaces}] Let $X$ be a simply connected symplectic
4-manifold and $\Sig$ a symplectically embedded surface of positive genus
and nonnegative self-intersection.
Assume also that $\pi_1(X\setminus\Sig)=1$. If
$K_1$ and $K_2$ are knots in $S^3$ and if
$(X, \Sig_{K_1})\cong (X,\Sig_{K_2})$, then  $\DD_{K_1}(t)=\DD_{K_2}(t)$.
Furthermore, if $\DD_K(t)\ne 1$, then
$\Sig_K$ is not smoothly ambient isotopic to a symplectic submanifold of $X$.
\end{thm}

The second part of the theorem holds because if $\Sig_K$ were symplectic,
$X\#_{\Sig_K=S_g}Y_g$ would be a symplectic manifold. The symplectic form
$\o$ on this
manifold is inherited from the forms on $X$ and $Y_g$; so $\la\o,R\ra=0$.
But $\sw_{X\#_{\Sig_K=S_g}Y_g}= \sw_{X\#_{\Sig=S_g}Y_g}\cdot\DD_K(r)$, and
it follows that the among the basic classes $k$ of $X\#_{\Sig_K=S_g}Y_g$,
more than one has $\la\o,k\ra$ maximal. This contradicts the fact that, for a
symplectic manifold, the maximality of $\la\o,K\ra$
characterizes the
canonical class among all basic classes \cite{T2}.

\section{Fiber sums of holomorphic Lefschetz fibrations \label{Lefschetz}}
In this section we shall construct for every integer $g \ge 3$ a pair
$(X_g, X'_g)$ of simply connected complex surfaces carrying holomorphic
genus $g$
Lefschetz fibrations with the property that their fiber sum (along a
regular fiber)
is a symplectic 4-manifold $Z_g$  which supports no complex structure;
in fact $Z_g$ is not even homeomorphic to a complex manifold.

Let $T(p,q)$ denote the $(p,q)$ torus knot in $S^3$ and let $N(p,q)$
denote the 4-manifold obtained by attaching a 2-handle to the 4-ball along
$T(p,q)$
with $0$-framing. It is well known that $N(p,q)$ is a Lefschetz fibration over
$D^2$ with generic fiber a Riemann surface of genus
$g(p,q)=(p-1)(q-1)/2$. Let
$W(p,q)$ denote the canonical resolution of the Brieskorn singularity
$\Sig(p,q,pq)$, the Seifert-fibered 3-manifold with three exceptional
fibers of order $p$, $q$, and $pq$, and with $H_1=\Z$. It is known that
$W(p,q)$ also supports the structure of a genus $g(p,q)$ Lefschetz
fibration
over $D^2$ with a singular fiber over $0$ which is a sequence of $2$-spheres
plumbed    according to the resolution diagram of $\Sig(p,q,pq)$. Finally, let
\[Z(p,q)= W(p,q) \cup N(p,q).\]
The manifold $Z(p,q)$ is a rational
surface  which is diffeomorphic to the connected sum of $\CP$ and $r(p,q)$
copies
of $\CPb$ for some computable integer $r(p,q)$. Furthermore, $Z(p,q)$
supports the
structure of a holomorphic Lefschetz fibration whose fiber has genus $g(p,q)$.

Now consider nontrivial torus knots $T(p,q)$ and $T(p',q')$ with the
property that
$g(p,q)=g(p',q')$. (This is  possible for every $g(p,q) \ge 3$.) Let
$F(p,q;p',q')$ denote the fiber sum along a regular fiber of
$Z(p,q)$ with $Z(p',q')$. Then $F(p,q;p',q')$ is a simply connected symplectic
$4$-manifold with
\[c_1^2=10+8g(p,q)-r(p,q)-r(p',q'), \;\;\chi=(b^++1)/2=1+g(p,q).\]

Furthermore, $F(p,q;p',q')$ supports the structure of a Lefschetz fibration
with fiber
of genus
$g(p,q)$. A computation of the Seiberg-Witten invariants of
$F(p,q;p',q')$  shows that, up to sign, there is a unique Seiberg-Witten basic
class. It follows that $F(p,q;p',q')$ is minimal.

\begin{conj} $F(p,q;p',q')$ supports the structure of a complex
4-manifold if and only if $\{p,q\}=\{p',q'\}$.
\end{conj}

\noindent As evidence, consider the pairs $(2, 2n+1)$ and $(3,n+1)$, $n
\not\equiv 2 \mod 3$. For
$F(2,2n+1;3,n)$ one can show that $r(2,2n+1)=4n+4$ and  $r(3,n+1)=3n+7$
so that
\[c_1^2=n-2, \;\;\chi=n+1. \] Thus, $c_1^2=\chi-3$, which violates the
Noether inequality $c_1^2 \ge 2\chi-6$. This means that
$F(2,2n+1;3,n)$ is a minimal symplectic 4-manifold that is not even
homotopy equivalent to a complex manifold.  In fact, it can be shown that
the fiber sum of $Z(2,2n+1)$ with itself is the elliptic surface
$E(n+1)$ and that the fiber sum of $Z(3,n+1)$ with itself is a Horikawa
surface with
$\chi = n+1$. Furthermore $F(2,2n+1;3,n)$ can be obtained from
$E(n+1)$ by removing from $Z(2,2n+1)\setminus F\subset E(n+1)$, $F$ a
regular fiber, the regular neighborhood of the configuration of  ($n-2$)
2-spheres:

\centerline{\unitlength 1cm
\begin{picture}(5,1.15)
\put(.9,.2){$\bullet$}
\put(1,.3){\line(1,0){1.3}}
\put(2.2,.2){$\bullet$}
\put(2.3,.3){\line(1,0){.75}}
\put(3.3,.3){.}
\put(3.5,.3){.}
\put(3.7,.3){.}
\put(4,.38){\line(1,0){.75}}
\put(4.65,.3){$\bullet$}
\put(.35,.6){$-(n+1)$}
\put(2.1,.6){$-2$}
\put(4.55,.6){$-2$}
\end{picture}}
\noindent whose boundary is the lens space $L((n-1)^2,-n)$ and replacing
it with the rational ball that this lens space  bounds. (See \cite{rat}
for all the details concerning this rational blowdown procedure.) Thus
$F(2,2n+1;3,n)$ is the manifold $Y(n)$ constructed in Lemma 7.5 of
\cite{rat}.

\section{Homeomorphic but non-diffeomorphic $4$-manifolds with the same
Seiberg-Witten invariants
\label{sameSW}}

In this section we construct examples of a pair $(X_1,X_2)$ of symplectic
4-manifolds with $X_1$ homeomorphic to $X_2$, $\sw_{X_1}=\sw_{X_2}$, but
$X_1$ is not
diffeomorphic to $X_2$. To do this choose a pair of fibered 2-bridge knots
$K(\alpha,\beta_1)$ and
$K(\alpha,\beta_2)$ with the same Alexander polynomials; for example
$K_1=K(105,64)$ and
$K_2=K(105,76)$ with Alexander polynomial
\[\Delta_K(t)=t^{-4}-5t^{-3}+13t^{-2}-21t^{-1}+25-21t+13t^2-5t^3+t^4.\]
Although these knots have the same Alexander polynomial, they can be
distinguished
by the fact that their branch covers are the lens spaces $L(\alpha,\beta_1)$
and
$L(\alpha,\beta_2)$ which are distinct; in our specific case $L(105,64)$ is
not diffeomorphic to $L(105,76)$. These knots are also
distinguished by their dihedral linking numbers; let
$S_{K_1}$ and
$S_{K_2}$ denote the 2-fold covers of
$S^3$ branched over
$K_1$ and $K_2$, with lifted branched loci $\tilde{K_1}$ and
$\tilde{K_2}$, respectively. Thus we have knots
$\tilde{K_i}$ in $S_{K_i}=L(\alpha,\beta_i)$. Take the $\alpha$-fold covers
of these
lens spaces to obtain links
$L_i=\{K_1^{(i)},\dots,K_\alpha^{(i)}\}$ which are
the lifts of the branch loci $\tilde{K_i}$. The linking numbers
of the links $L_1$ and $L_2$ are known as the `dihedral linking numbers' of
the 2-bridge knot $K(\a,\b)$.

Now perform the knot surgery construction of \S\ref{knots} on the $K3$
surface, replacing $T^2\times D^2$ with
$S^1\times(S_{K_j}\setminus\tilde{K_j})$. The resulting 4-manifolds are
the manifolds $X_i$. Either by adapting the arguments of
\cite{KL4M} or by using \cite{IP} or \cite{L}, it can be shown that
$\sw_{X}=\sw_{Y}=\Delta_K(t)\cdot \Delta_K(-t)$.
Unfortunately, the $X_i$ are not simply connected (but are
homeomorphic). In particular,
$\pi_1(X_1)=\pi_1(X_2)=\Z_{\a}$, and the
$\alpha$-fold covers $\tilde{X_1}$ and $\tilde{X_2}$ of $X_1$ and $X_2$ are not
diffeomorphic. To see this, observe that $\tilde{X_i}$ is
obtained as our link construction in \cite{KL4M} (cf. \S~\ref{knots}) by
gluing one
copy of
$E(2)$ minus a neighborhood of a smooth elliptic fiber to every boundary
component of
$S^1\times (S^3\setminus L_i)$. It follows from \cite{KL4M} that
\[ \ssw_{\tilde{X_i}}=\Delta_{L_i}(t_1,\dots,t_\alpha)\cdot
    \prod_{j=1}^{\a} (t_j^{1/2}-t_j^{-1/2})\]
\noindent Since the linking numbers of the links $L_1$ and $L_2$ are
different, it can be shown the Hosokawa polynomials of the links
$L_1$ and $L_2$, when evaluated at $1$ are distinct \cite{FSHos}. Thus their
Alexander polynomials are different and $\tilde{X_1}$ is not diffeomorphic to
$\tilde{X_2}$.

There is a lesson to be learned from these examples. One must consider the
Seiberg-Witten invariants of a 4-manifold $X$ together with those of all of
its covers as the appropriate invariant for $X$.

\section{Nonsymplectic 4-manifolds with one basic class}

Recall from \S~\ref{knots}, that if $k$ is a basic class of $X$, so is $-k$.
Because of this, we say that $X$ {\it{has $n$ basic classes}} if the set $\{k\,
|\,\ssw_X(k)\ne0\}/ \{\pm1\}$ consists of $n$ elements.
There are abundant examples of $4$-manifolds with one basic class. Minimal
nonsingular
algebraic surfaces of general type have one basic class (the canonical class)
\cite{W}. The authors and others have constructed many examples of minimal
symplectic
manifolds with one basic class and $\chi -3\le {c_1}^2 < 2\chi -6$. (These
manifolds
cannot admit complex structures due to the geography of complex surfaces.)
However, the
examples described here are the first nonsymplectic manifolds with one basic
class.

Let $X=E(2)$ and $T$ a smooth elliptic fiber. For a knot $K$ of genus $g$
form the knot
surgery construction to obtain $X_K$. In $X_K$ there is a surface
$\Sig$ of genus $g+1$ with $[\Sig]^2=0$ and $[\Sig]\cdot [T]=1$. Let $M$ be the
3-manifold obtained from 0-surgery on the trefoil knot. Then $S^1\x M$ is a
$T^2$-fiber
bundle over $T^2$. The fiber sum of $g+1$ copies of the fiber bundle gives
a 4-manifold
$Y$ which is an $F=T^2$-bundle over a surface of genus $g+1$, and it is
easily seen that
there is a section $C$. Furthermore, $Y$ is a symplectic 4-manifold
with $c_1(Y)=-2g[F]$. Our example, corresponding to the genus $g$ knot $K$ is
$ Z_K=X_K\#_{\Sig=C}Y$.
We perform this fiber sum so that $Z_K$ is a spin 4-manifold \cite{Gompf}.
It can be
seen to be simply connected.

Write the symmetrized Alexander polynomial of $K$ as
$\DD_K(t)=a_0+\sum\limits_{n=1}^da_n(t^n+t^{-n})$, and call $d$ the
{\it{degree}} of
$\DD_K(t)$. Since the genus of $K$ is $g$, we have
$d\le g$. If $K$ is an alternating knot, for example, then $d=g$. Say that the
Alexander polynomial of $K$ has {\it{maximal degree}} if $d=g$. Using
techniques of
\cite{MST} we calculate:

\begin{thm}[\cite{1basic}] Let $K$ be a knot in $S^3$ whose Alexander
polynomial has
maximal degree. Then $Z_K$ has one basic class, $k$, with
$|\ssw_{Z_K}(k)|=a_d$, the top coefficient of
$\DD_K(t)$. When $|a_d|>1$, $Z_K$ is nonsymplectic.
\end{thm}

%%--------------------Here the manuscript ends--------------------------------

\Addresses
\end{document}